\documentstyle[12pt]{amsart}
\newtheorem{theorem}{Theorem}[section]

\newtheorem{lemma}[theorem]{Lemma}

\newtheorem{claim}{Claim}[lemma]
\newtheorem{problem}{Problem}

\theoremstyle{definition}

\newtheorem{definition}[theorem]{Definition}

\theoremstyle{remark}

\newcommand{\proof}{\begin{pf}}
\newcommand{\Proof}[1]{\begin{pf*}{Proof of #1}}
\newcommand{\eproof}{\end{pf}}
\newcommand{\Eproof}{\end{pf*}}

\newcommand{\arablabel}{
          \renewcommand{\labelenumi}{{\rm (\arabic{enumi})}}
          \renewcommand{\theenumi}{{\rm (\arabic{enumi})}}
                    }

\newcommand{\alabel}{
          \renewcommand{\labelenumi}{{\rm (\alph{enumi})}}
          \renewcommand{\theenumi}{{\rm (\alph{enumi})}}
                    }
\newcommand{\Alabel}{
          \renewcommand{\labelenumi}{{\rm (\Alph{enumi})}}
          \renewcommand{\theenumi}{{\rm (\Alph{enumi})}}
                    }
\newcommand{\rlabel}{
          \renewcommand{\labelenumi}{{\rm (\roman{enumi})}}
          \renewcommand{\theenumi}{{\rm (\roman{enumi})}}
                    }

\newcommand{\nolabel}{
         \renewcommand{\labelenumi}{}
         \renewcommand{\theenumi}{}}

\newcommand{\witem}{\renewcommand{\labelitemi}{{\bf (W)}}}

\newcommand{\bcal}{{\cal B}}

\newcommand{\gcal}{{\cal G}}
\newcommand{\hcal}{{\cal H}}
\newcommand{\ical}{{\cal I}}

\newcommand{\ncal}{{\cal N}}
\newcommand{\pcal}{{\cal P}}

\newcommand{\ucal}{{\cal U}}

\newcommand{\wcal}{{\cal W}}

\newcommand{\setm}{\setminus}
\newcommand{\empt}{\emptyset}
\newcommand{\subs}{\subset}

\newcommand{\oo}{{{\omega}_1}}
\newcommand{\rest}{\lceil}
\newcommand{\dom}{\mbox{\rm{dom}}}

\def\<{\left\langle}
\def\>{\right\rangle}

\def\OO{{\omega}}
\def\oo{\omega_1}
\def\br#1;#2;{\bigl[ {#1} \bigr]^ {#2} }
\def\bc#1;#2;{\bigl( {#1} \bigr)^ {#2} }

\def\ooseq#1;#2;{\< {#1}_{#2}:{#2}<\oo\>}
\def\ooset#1;#2;{\{ {#1}_{#2}:{#2}<\oo\}}
\def\seq#1;#2;#3;{\< {#1}_{#2}:{#2}<#3\>}
\def\set#1;#2;#3;{\{ {#1}_{#2}:{#2}<#3\}}
\def\oseq#1;#2;{\< {#1}_{#2}:{#2}<\OO\>}
\def\oset#1;#2;{\{ {#1}_{#2}:{#2}<\OO\}}
\def\oosequ#1;#2;{\< {#1}^{#2}:{#2}<\oo\>}
\def\oosetu#1;#2;{\{ {#1}^{#2}:{#2}<\oo\}}
\def\sequ#1;#2;#3;{\< {#1}^{#2}:{#2}<#3\>}
\def\setu#1;#2;#3;{\{ {#1}^{#2}:{#2}<#3\}}
\def\osequ#1;#2;{\< {#1}^{#2}:{#2}<\OO\>}
\def\osetu#1;#2;{\{ {#1}^{#2}:{#2}<\OO\}}

\def\force{\raisebox{1.5pt}{\mbox{$\scriptscriptstyle\|$}}
\mbox{$\!\mbox{---}$}}

\def\to{\longrightarrow}

\def\w{\operatorname{w}}
\def\nw{\operatorname{nw}}
\def\RR{\operatorname{R}}
\def\good{good }

\def\goodi#1;{$#1$-\good}
\def\nea{neighbourhood assignment }
\def\dgf{D_G^f}
\def\dnf{D_N^f}
\def\dhf{D_H^f}
\def\ugx{U(G,x)}
\def\uhy{U(H,y)}

\def\anfi#1;{\<A^{#1},n^{#1},f^{#1},g^{#1},k^{#1},d^{#1},e^{#1}\>}
\newcommand{\anfg}{\<A,n,f,g\>}
\def\anfgi#1;{\<A^{#1},n^{#1},f^{#1},g^{#1}\>}
\newcommand{\bde}{\<B,d,e\>}
\def\bdei#1;{\<B^{#1},d^{#1},e^{#1}\>}
\newcommand{\bdne}{\<B,m,d,e\>}
\def\bdnei#1;{\<B^{#1},d^{#1},m^{#1},e^{#1}\>}

\def\ap{A^p}
\def\fp{f^p}

\def\np{n^p}

\def\up{U^p}
\def\sbar{\overline{\sigma}}
\def\sstr{{\sigma}^*}

\newcommand{\id}{\operatorname{id}}

\newcommand{\tle}{\triangleleft}
\newcommand{\ai}{{\alpha},i}
\newcommand{\bj}{{\beta},j}
\newcommand{\cl}{{\gamma},l}
\newcommand{\sai}{\sbar({\alpha}),i}
\newcommand{\sbj}{\sbar({\beta}),j}
\newcommand{\scl}{\sbar({\gamma}),l}

\newcommand{\sba}{\sbar({\alpha})}
\newcommand{\sbb}{\sbar({\beta})}
\newcommand{\sbc}{\sbar({\gamma})}
\newcommand{\ssa}{\sstr({\alpha})}

\newcommand{\ssc}{\sstr({\gamma})}

\newcommand{\Pko}{P^{\kappa}}
\newcommand{\Pkoo}{P^{\kappa}_0}
\newcommand{\pcalk}{\pcal^{\kappa}}

\newcommand{\qsd}{Q_0}
\newcommand{\qssc}{Q_1}

\newcommand{\cont}{2^{\omega}}
\newcommand{\irr}{\<\ncal,{\varepsilon}\>}

\def\fgi#1;{\<f_#1,g_#1\>}

\def\fghi#1;{\<f_{#1},g_{#1},h_{#1}\>}

\newcommand{\dd}{\operatorname{d}}
\newcommand{\hh}{\operatorname{h}}
\def\dyi#1;{D[#1]}
\newcommand{\dy}{\dyi y;}
\theoremstyle{definition}
\newtheorem{ocaaxiom}{Open Coloring Axiom}

\title{What makes a space have large weight?}
\thanks{The preparation of this paper was supported by the 
Hungarian National Foundation for Scientific Research grant no. l908 }

\author{I. Juh\'asz}
\address{Mathematical Institute of the Hungarian Academy of Sciences}
\email{h1152juh@@ella.hu}
\author{L. Soukup}
\thanks{The second author  was supported by the 
Magyar Tudom\'any\'ert Foundation}
\address{Mathematical Institute of the Hungarian Academy of Sciences}
\email{h1153sou@@ella.hu}
\author{Z. Szentmikl{\'o}ssy}
\address{E{\"o}tv{\"o}s Lor{\'a}nd University, Department of Analysis}

\begin{document}
%\begin{comment}
\maketitle
\begin{abstract}
In section \ref{sc:irred} of this paper we formulate several conditions
(two of them are necessary and sufficient) which imply that a space 
of small character has large weight.
In section \ref{sc:example} we construct a ZFC example of a 0-dimensional space
$X$ of size $\cont $ with $\w(X)=\cont$ and ${\chi}(X)=\nw(X)={\omega}$,
 we show that CH implies the existence of a 0-dimensional space $Y$
of size $\oo$ with $\w(Y)=\nw(Y)=\oo$ and ${\chi}(Y)=\RR(Y)={\omega}$,
and we prove that it is consistent that $\cont$ is as large as you wish
and there is a 0-dimensional space $Z$ of size $\cont$ such that  
$\w(Z)=\nw(Z)=\cont$ but  ${\chi}(Z)=\RR(Z^{\omega})={\omega}$.
\end{abstract}

\section{Introduction}

Since ${\chi}(X)\ge |X|$ implies $\w(X)={\chi}(X)$, one possible answer
to the question in the title is that having large character will make a space
have large weight. 
Thus we arrive at the following more interesting problem:
What makes a space have weight larger than its character?
Discrete spaces give examples of such spaces but the Sorgenfrey line
is first countable, has weight $\cont $ but it has no uncountable discrete 
subspace. The reason for the latter space to have weight $\cont$
is that it is weakly separated, i.e.,
one can assign to every point $x$ a neighbourhood $U_x$ such that 
$x\ne y$ implies either $x\notin U_y$ or $y\notin U_x$. So we may ask now whether every first countable space of ``large'' weight has a ``large''
weakly separated subspace? This question was the actual starting point
of our investigations, and while we found a negative answer to it we also
succeeded in finding successively more and more general conditions that ensure  
having large weight for spaces of small character.

In section \ref{sc:irred} we introduce the notion of 
{\em irreducible base of a space} (see definition
\ref{df:irred}) and investigate its basic properties.
This notion is a weakening of weakly separatedness but
the existence of such a base still implies that the 
weight of the space can not be smaller than its cardinality.
The main advantage of this notion, in contrast to weakly separatedness,
  lies in the fact that, 
as we will see in section \ref{sc:example},
a large space with an irreducible
base might have small net weight.

This leads to the formulation of the following problem:
\begin{problem}
Does every  first countable space of uncountable weight 
contain an uncountable subspace with an irreducible base?
\end{problem}

In section \ref{sc:example} we construct examples. 
First a ZFC example is given of a space $Y$
with $|Y|=\w(Y)=\cont$ and ${\chi}(Y)=\RR(Y)=\nw(Y)={\omega}$.
After seeing that ${\chi}(Y)\RR(Y)<\w(Y)$  but 
${\chi}(Y)\RR(Y)\ge\nw(Y)$    in the above mentioned 
example, we asked whether $\nw(X)\le\RR(X){\chi}(X)$
or just $\nw(X)\le\RR(X^{\omega}){\chi}(X)$ are provable 
for every $T_2$ or  regular space $X$. 
Using CH a 0-dimensional counterexample is given  
to the first question and   using a c.c.c forcing argument 
we disprove the second inequality in section \ref{sc:example}. 
However we don't know ZFC counterexamples.  
\begin{problem}
Is there  a ZFC example of a space $X$ satisfying
 $\RR(X^{\omega})={\chi}(X)={\omega}$ but $\nw(X)>{\omega}$?
\end{problem}
We know that under ${MA}$ the cardinality of such a space
must be at least $\cont$ (see \cite{JSS2}).
In \cite[p 30]{To} Todor{\v c}evi{\v c} introduced the axiom (W):
\begin{itemize}
\witem
\item {\em If $X$ is  regular space with $\RR(X)^{\omega}={\omega}$ then
$\nw(X)={\omega}$}.
\end{itemize}
and he claimed that PFA implies (W).

We use  standard topological notation and terminology throughout, cf \cite{J}.

\section{Conditions ensuring large weight}
\label{sc:irred}
\begin{definition}
Given a topological space $\<X,{\tau}\>$ and a subspace $Y\subs X$ a function $f$
is called a {\em \nea on Y} iff $f:Y\to {\tau}$ and $y\in f(y)$ for each 
$y\in Y$. 
\end{definition}

The notion of weakly separated spaces and the cardinal function
$\RR$ were introduced by Tka{\v c}enko in \cite{Tk}.

\begin{definition}
A space $Y$ is {\em weakly separated} if we can find
a \nea $f$ on $Y$ such that
$$
\forall y\ne z\in Y\ ( y\notin f(z) \lor z\notin f(y) ),
$$
moreover
$$
\RR(X)=\sup\{|Y|:Y\subs X\mbox{ is weakly separated}\}.
$$
\end{definition}

Obviously $\RR(X)\le \nw(X)$.
Tka{\v c}enko asked whether $\RR(X)=\nw(X)$ is provable for regular spaces.
Hajnal and Juh{\'a}sz, in \cite{HJ}, gave several 
consistent counterexamples
using CH and some c.c.c forcing arguments.
However, their spaces were not first countable.

If one wants to construct a  first countable space 
 on $\oo$ without uncountable weakly separated subspaces 
a natural idea is to force with  finite approximations of 
a  base of such a space.  The space  $X$ given by a generic
filter satisfies $\RR(X)={\omega}$, but without additional assumptions
standard density arguments give $\w(X)={\omega}$, too. 
To ensure large weight of the generic space we 
actually needed that the base should satisfy a certain property.
As it turned out this notion proved  to be  useful not only 
in the special forcing construction.
Its definition is now given below.

\begin{definition}\label{df:irred}
Let $X$ be a topological space. 
A base $\ucal$ of $X$ is called {\em irreducible} 
if it has an {\em irreducible decomposition} $\ucal=\bigcup\{\ucal_x:x\in X\}$,
i.e,  $(i)$ and $(ii)$ below hold:
\begin{enumerate}\rlabel
\item $\ucal_x$ is a neighbourhood base of $x$ in $X$ for each $x\in X$.
\item for each $x\in X$   the family 
$\ucal^-_x=\bigcup\limits_{y\ne x}\ucal_y$ is not a base of $X$,
hence it doesn't contain a neighbourhood base of 
$x$ in $X$. 
\end{enumerate}
\end{definition}

Let $\ucal$ be an irreducible base with the irreducible 
decomposition $\{\ucal_x:x\in X\}$. Then for each $x\in X$,
since  $\bigcup\limits_{y\ne x}\ucal_y$ does not contain a neighbourhood base of 
$x$ in $X$,  we can fix an open neighbourhood $U_x$ such that 
$$
\forall y\in X\setm \{x\}\ \forall V\in \ucal_y
\ ( x\in V\ \longrightarrow \ V\setm U_x\ne\empt).
$$
Let $\ucal^*_x=\{U\in\ucal_x:U\subs U_x\}$. Then 
$\ucal^*=\bigcup\{\ucal^*_x:x\in X\}$ is an irreducible base of $X$ and
its irreducible decomposition $\{\ucal^*_x:x\in X\}$ has  the following
 property $(*)$:
\begin{equation}\tag{$*$}
\forall x\ne y\in X\ \forall U\in \ucal^*_x
\ \forall V\in\ucal^*_y\ ( x\in V\land y\in U )\longrightarrow  V\setm U\ne\empt.
\end{equation}

To simplify our notation we will say that  a base 
$\ucal$ has property $(*)$ if it has a decomposition 
$\ucal=\bigcup\{\ucal_x:x\in X\}$ satisfying $(i)$ and $(*)$ above.
Obviously,  any base with  property $(*)$ is irreducible. 
So we established the following lemma:
\begin{lemma}
A space $X$ has an irreducible base iff it has a base with property
$(*)$.
\end{lemma}
 
The next two lemmas establish the basic connection between 
weakly separatedness, existence of irreducible base and 
the requirement $\w(X)\ge |X|$.
\begin{lemma}
If $X$ is weakly separated, then $X$ has an irreducible base.
\end{lemma}

\proof
Let $f$ be a \nea on $X$ witnessing that it is weakly separated.
Take $\ucal_x=\{G\in{\tau}_X:x\in G\subs f(x)\}$ and 
$\ucal=\bigcup\{\ucal_x:x\in X\}$. Then 
$\{\ucal_x:x\in X\}$ is an irreducible decomposition of the base
$\ucal$.
\eproof

The converse of this lemma fails as we will see it later 
(theorem \ref{tm:zfc_example}).

\begin{lemma}
If $X$ has an irreducible base, then $\w(X)={\chi}(X)|X|$.
\end{lemma}

\proof
If ${\chi}(X)\ge |X|$ this is trivial, so assume that
${\lambda}={\chi}(X)<|X|$.
Consider an irreducible base $\ucal$ with irreducible decomposition 
$\{\ucal_x:x\in X\}$. We can assume that 
$|\ucal|=\w(X)$ and $|\ucal_x|\le{\lambda}$ for each $x\in X$.
If $\wcal\subs \ucal$ with $|\wcal|<|X|$, then there is
$x\in X$ with $\wcal\cap\ucal_x=\empt$, so $\wcal$ can't be
a base by the irreducibility of $\ucal$. Thus
$\w(X)=|\ucal|\ge |X|$. 
\eproof

\begin{definition}
Given a topological space $X$, a subspace $Y\subs X$,
a \nea $f$ on $Y$ and a set $N\subs X$ let
$$
\dnf=\{y\in Y: y\in N\subs f(y)\}.
$$ 
\end{definition}

The following results show that both weakly separatedness and having
an irreducible base may be characterized with the existence of a \nea 
$f$ such that $\dgf$ is ``small'' in some sense for each open $G$.
For example we have the following easy result  whose proof we leave to the
 reader.

\begin{theorem}
\label{tm:dgf_wseparated}
Given a topological space $X$, a subspace $Y\subs X$ is weakly separated
iff there is a \nea $f$ on $Y$ such that $|\dgf|\le 1$ for each
open $G\subs X$. \qed
\end{theorem}

\begin{lemma}
\label{lm:dgf_discrete}
If a space  $X$ has an irreducible base, then 
there is a \nea  $f$ on $X$ such that 
$\dgf$ is closed and discrete in $G$ for all open $G\subs X$.
\end{lemma}

\proof
Let $\ucal$ be a  base of $X$ having a decomposition $\{\ucal_x:x\in X\}$ 
with property $(*)$ and fix a \nea  $f$ with $f(x)\in \ucal_x$.
Assume on the contrary that $x\in G$ is an accumulation point of
$\dgf$ for some open $G\subs X$.  Choose $U\in\ucal_x$ with $x\in U\subs G$.
Pick $y\in \dgf\cap U$, $y\ne x$. Then $y\in U\in \ucal_x$,
$x\in G\subs f(y)\in\ucal_y$, and $U\subs f(y)$, which contradicts 
property $(*)$ of $\ucal$.
\eproof

\begin{theorem}
\label{tm:irred_dgf} The following statements are equivalent for any regular
space $\<X,{\tau}\>$:
\begin{enumerate}\arablabel
\item $X$ has an irreducible base.
\item There is a \nea  $f$ on $X$ such that 
$\dgf$ is closed and discrete in $G$ for all open $G\subs X$.
\item There is a \nea  $f$ on $X$ such that 
$\dgf$ is a discrete subspace of $X$ for each open $G\subs X$.
\end{enumerate}
\end{theorem}

\proof
\mbox{}\newline
{\bf $(1)\to (2)$.}
\newline 
This is just lemma \ref{lm:dgf_discrete}.
\newline 
{\bf $(2)\to (3)$.}\newline 
 Straightforward.
\newline 
{\bf $(3)\to (1)$.}\newline 
 Fix a \nea on $X$ witnessing $(3)$. 
Since $X$ is regular we can assume that $f(x)$ is regular open for each 
$x\in X$ and that $f(x)=\{x\}$ provided  $x$ is isolated.
Given an open $G\subs X$ set
$$
\ugx=(G\setm\overline{\dgf})\cup\{x\}
$$
and put 
$$
\ucal_x=\{\ugx:x\in \dgf \land\mbox{$G$ is regular open}\}.
$$
Since $\dgf$ is discrete and $x\in\dgf$, we have that $\ugx$ is also open, 
and $\ucal_x$ is a neighbourhood base of $x$ because
$x\in\ugx\subs G$.
We claim that $\ucal=\bigcup\{\ucal_x:x\in X\}$ is an irreducible base
because the decomposition $\{\ucal_x:x\in X\}$ has property $(*)$.
Assume on the contrary that
  $x\ne y\in X$, $\ugx\in \ucal_x$,
$\uhy\in \ucal_y$ with  $\{x,y\}\subs \ugx\cap\uhy$ and
$\ugx\subs \uhy$.
Since $|G|>1$ and $|f(x)|=1$ whenever $x$ is isolated in $X$, it
follows that  $\dgf$  can not contain isolated points from $X$.
But this set is also discrete, so  $\overline{\dgf}$ is nowhere dense in $X$.
Since  $H$ is regular open,
$\ugx=(G\setm\overline{\dgf})\cup\{x\}\subs H$ implies $G\subs H$.
Thus $y\in \dgf$ for $y\in\ugx\subs G\subs H\subs f(y)$ and so $y\notin \ugx$,
which is impossible.
\eproof

We don't know if the assumption on the regularity of $X$
is essential in theorem \ref{tm:irred_dgf}.

Next we show that the existence of an $f$
with $\dgf$ ``small'' for all open sets $G$
already implies that the weight of our space is large.

\begin{definition}
A topological space $Y$ is {\em pseudo weakly separated }
if it contains a weakly separated subspace $Z$ with $|Z|=|Y|$.
\end{definition}

\begin{theorem}
\label{tm:c_aa}
Let  $X$ be a topological space, $Y\subs X$,
 $f$ be  a \nea  on   $Y$ 
and ${\lambda}\le |X|$ be a regular cardinal. 
If  $\dgf$ is  the union of $<{\lambda}$ many pseudo weakly separated
subspaces for each open $G\subs X$, then $\w(X)\ge |Y|$.  
\end{theorem}

\proof
Assume on the contrary that $\bcal$ is a base with 
$|\bcal|<|Y|$ and let ${\kappa}=|\bcal|^++{\lambda}$.
 Since $Y=\bigcup\limits_{G\in\bcal}\dgf$,
there is a $G\in\bcal$ with $|\dgf|\ge{\kappa}$. 
But $\dgf$ is the union of $<{\lambda}$ many pseudo weakly separated
subspaces, so one of them has cardinality $\ge {\kappa}$.
Thus $\dgf$ contains a weakly separated subspace $Z$ with 
$|Z|\ge{\kappa}$. Hence 
$\w(X)\ge \w(Y)\ge \w(Z)\ge {\kappa}>|\bcal|\ge\w(X)$, 
which is impossible.
\eproof

Since weakly separated spaces have not just large weight but also large
net weight,
 if we assume that $\dnf$ is like in theorem \ref{tm:c_aa}
for all subsets $N\subs X$, the same argument yields that even 
the net weight of $X$ is large.

\begin{theorem}
\label{tm:c_aaa}
Let $f$ be a \nea on a topological space $X$ and ${\lambda}\le |X|$ be
a regular cardinal. 
If  $\dnf$ is  the union of $<{\lambda}$ many pseudo weakly separated
subspaces for each $N\subs X$,  then $\nw(X)\ge |X|$.
\qed
\end{theorem}

The following results show that the above type of 
``smallness'' assumptions on $\dgf$ can  actually be used
to characterize spaces of small character and large weight!

\begin{theorem}\label{tm:wk1}
Let  ${\kappa}$ be a cardinal and $X$ a topological space with  
${\chi}(X)<{\kappa}$. 
(If ${\kappa}$ is regular then the assumption ${\chi}(p,X)<{\kappa}$ 
for each $p\in X$ would suffice.)
 Then the following are equivalent:
\begin{enumerate}\alabel
\item $\w(X)\ge{\kappa}$,
\item There is a subspace $Y\subs X$ of size ${\kappa}$ and a \nea
$f$ on $Y$ such that $\dgf$ is right-separated for each open 
$G\subs X$.
\item There is a subspace $Y\subs X$ of size ${\kappa}$, a \nea
$f$ on $Y$ and a regular cardinal ${\lambda}\le {\kappa}$ 
such that $\dgf$ is  the union of $<{\lambda}$ many pseudo weakly separated 
subspaces for each open $G\subs X$.
\end{enumerate}
\end{theorem}

\proof
{\mbox{}}\newline 
{\bf (a)$\to$(b).}\newline 
 For each $x\in X$ fix a neighbourhood base $\bcal_x$
of $x$ in $X$ with minimal cardinality. Since 
$\w(X)\ge {\kappa}$ we can  construct a sequence 
$\{y_{\eta}:{\eta}<{\kappa}\}\subs X$ such that
for each ${\eta}<{\kappa}$ the family
 $\bigcup\limits_{{\xi}<{\eta}}\bcal_{y_{\xi}}$ does not contain a base 
of $y_{\eta}$ in $X$ and we can pick an open set 
$f({\eta})\in {\tau}_X$ which witnesses
this, i.e., $y_{\eta}\in f({\eta})$ and there is no
$U\in\bigcup\limits_{{\xi}<{\eta}}\bcal_{y_{\xi}}$ with
$y_{\eta}\in U\subs f({\eta})$. We claim that the \nea  $f$ on
$Y=\{y_{\eta}:{\eta}<{\kappa}\}$ has the property that
$\dgf$ is right separated in its natural order for each open $G$.
Assume on the contrary that there is an open $G$ and ${\xi}<{\kappa}$
such that $y_{\xi}\in\dgf\cap\overline{\{y_{\eta}\in\dgf:{\eta}>{\xi}\}}$.
Since  $H\subs G$ implies  $\dgf\cap H\subs \dhf$ we can assume that 
$G\in\bcal_{y_{\xi}}$. Then  there is ${\eta}>{\xi}$
with $y_{\eta}\in G\cap \dgf$. Hence $y_{\eta}\in G\subs f({\eta})$
and $G\in\bcal_{y_{\xi}}$,  contradicting  the choice of 
$y_{\eta}$ and  $f({\eta})$. 
\newline 
{\bf (b)$\to$(c).}
\newline   Straightforward.
\newline   
{\bf (c)$\to$(a)} 
\newline This is immediate from  theorem \ref{tm:c_aa}.
\eproof

Conditions (b) and (c) in \ref{tm:wk1} have the
(perhaps just aesthetic) drawback that the requirements on the subspace $Y$
are external in nature, i.e. they do not only depend on $Y$.
This drawback is eliminated in the following result, which however works only 
for regular spaces and regular cardinals.

\begin{theorem}\label{tm:wk2}
Let $X$ be a regular topological space and ${\kappa}$  a regular
cardinal with ${\chi}(p,X)<{\kappa}$ for each $p\in X$. 
Then the following are equivalent:
\begin{enumerate}\alabel
\item $\w(X)\ge{\kappa}$,
\item There is a subspace $Y\subs X$ of size ${\kappa}$ and a \nea
$f:Y\to {\tau}_Y$ such that $|\dgf|<{\kappa}$  for each open 
$G\in {\tau}_Y$.
\item There is a subspace $Y\subs X$ of size ${\kappa}$ and a \nea
$f:Y\to {\tau}_Y$ such that $\dgf$  is the union of
$<{\kappa}$ many pseudo weakly separated subspaces for each open 
$G\in {\tau}_Y$.
\end{enumerate}
\end{theorem}

\proof
\mbox{}\newline 
{\bf (a)$\to$(b).}\newline 
If there is a weakly separated subspace $Y\subs X$ with $|Y|={\kappa}$ 
then we are done by theorem \ref{tm:dgf_wseparated}. Otherwise we have
$\hh(X)<{\kappa}$ and  
$\dd(X)<{\kappa}$, hence we can pick a dense $D\subs X$ with $|D|<{\kappa}$.
For each $x\in X$ fix a neighbourhood base $\bcal_x$
of $x$ in $X$ with minimal cardinality containing only regular open sets. 
Since $\w(X)\ge {\kappa}$ we can  construct a sequence 
$\{y_{\eta}:{\eta}<{\kappa}\}\subs X\setm D$ such that
for each ${\eta}<{\kappa}$ writing 
$D_{\eta}=\{y_{\xi}:{\xi}<{\eta}\}$ the family
 $\bigcup\{\bcal_y:y\in D\cup D_{\eta}\}$ 
does not contain  a base 
of $y_{\eta}$ in $X$ and we can pick a regular open set 
$f'({\eta})$ which witnesses this, i.e., $y_{\eta}\in f'({\eta})$ and there is no $U\in \bigcup\{\bcal_y:y\in D\cup D_{\eta}\}$ with 
$y_{\eta}\in U\subs f'({\eta})$. 
Let $Y=D\cup \{y_{\eta}:{\eta}<{\kappa}\}$ and define the \nea $f$ 
on $Y$ by $f(y_{\eta})=f'(y_{\eta})\cap Y$ for ${\eta}<{\kappa}$
and $f(d)=Y$ for $d\in D$.
We claim that the space
$\<Y,{\tau}_Y\>$ and the \nea  $f$ have the property that
$Z_G=\dgf\setm D$ is right separated in the  order inherited from the
indexing for each open $G$.
Assume on the contrary that there is an open $G$ and ${\xi}<{\kappa}$
such that $y_{\xi}\in Z_G\cap\overline{\{y_{\eta}\in Z_G:{\eta}>{\xi}\}}$.
Since  $H\subs G$ implies  $\dgf\cap H\subs \dhf$ we can assume that 
there is $G'\in\bcal_{y_{\xi}}$ with $G=G'\cap Y$. 
Then  there is ${\eta}>{\xi}$
with $y_{\eta}\in G\cap Z_G$. 
But $D\subs Y$ is dense in $X$ and $G'$ is regular open, so
 $y_{\eta}\in G\subs f({\eta})$ implies 
$G'\subs f'({\eta})$ which contradicts  the choice of 
$y_{\eta}$ and  $f'({\eta})$. 
Thus $Z_G$   is right separated and so $|Z_G|<{\kappa}$.
Therefore $|\dgf|=|Z_G|+|D|<{\kappa}$, which means that 
$\<Y,{\tau}_Y\>$ and the \nea  $f$ satisfy (b).
\newline 
{\bf (b)$\to$(c).}\newline 
Straightforward. \newline
{\bf (c)$\to$(a).}\newline 
$\w(X)\ge\w(Y)\ge{\kappa}$ by theorem \ref{tm:c_aa}. 
\eproof

\section{Examples of spaces with large weight and small character}
\label{sc:example}

Denote by $\irr$ the space of irrational numbers endowed with the Euclidean
topology.  For $x\in\ncal$ and ${\eta}>0$ write
$U(x,{\eta})=(x-{\eta},x+{\eta})\cap \ncal$.

\begin{theorem}
\label{tm:zfc_example}
There is a set $X\subs \ncal$ of size $\cont$ and a 0-dimensional 
first countable refinement ${\tau}$ of ${\varepsilon}$ on $X$ such that
\begin{enumerate}\rlabel
\item $X=Y\cup Z$, where ${\tau}_Y={\varepsilon}_Y$ and
${\tau}_Z={\varepsilon}_Z$,
\item $\<X,{\tau}\>$ has an irreducible base.
\end{enumerate}
Thus ${\chi}(X)=\nw(X)={\omega}$ but $\w(X)=\cont$.
\end{theorem}

\proof
Let $Z=\{z_n:n\in {\omega}\}\subs \ncal$ be dense.
Fix a nowhere dense closed set $Y\subs\ncal\setm Z$
of size $\cont$. Let $X=Y\cup Z$.
For each $y\in Y$ choose a  strictly increasing sequence 
of pairwise disjoint
intervals with rational endpoints, $\ical^y=\{I^y_n:n\in {\omega}\}$,
such that $\ical^y$ converges to $y$ and   $J^y=\bigcup\ical^y$ 
is disjoint from $Y$.
This can be done because $Y$ is nowhere dense.
Set $J^z=\empt$ for $z\in Z$.
For $x\in Y$ and ${\eta}>0$ let
  $V(x,{\eta})=(U(x,{\eta})\setm  J^x)\cap X$.
Let the neighbourhood base of $x\in Y$ in ${\tau}$ be 
$$
\bcal_x=\{V(x,{\eta}):{\eta}>0\}.
$$
If $z=z_n\in Z$ then pick ${\eta}_n>0$ such that 
$U(z_n,{\eta}_n)$ is disjoint from $Y\cup \{z_i:i<n\}$
and put 
$$
\bcal_{z_n}=\{U(x,{\eta}):{\eta}_n>{\eta}>0\}.
$$
Since 
\begin{equation}
\tag{\dag}\forall {\eta}>0\ \forall x\in V(y,{\eta})\setm \{y\}
\ \exists {\delta}>0\ U(x,{\delta})\cap X\subs V(y,{\eta})
\end{equation}
it follows that $\bcal=\bigcup\{\bcal_x:x\in X\}$ is a base of a topology.
We claim that $\bigcup\{\bcal_x:x\in X\}$ is an irreducible decomposition
of $\bcal$ because it has property $(*)$. 
So let $u,v\in X$, $U\in \bcal_u$, $V\in\bcal_v$ with 
$\{u,v\}\subs U\cap V$. Then $u,v\in Y$ because
$W\in \bcal_{z_n}$ implies $W\cap X\subs \{z_k:k\ge n\}$.
The density of $Z$ in $\ncal$  implies
\begin{equation}
\tag{+}\forall {\eta}>0\ \forall{\delta}>0\ 
U(y,{\delta})\cap X\not\subs V(y,{\eta})
\end{equation}
for each $y\in Y$. 
But $v\in U$ implies that there is some 
${\eta}$ with $U(v,{\eta})\subs U$, 
so $U\setm V\ne\empt$. Thus the 
base $\bcal$ is  irreducible.
On the other hand, ${\tau}_Y={\varepsilon}_Y$ 
and ${\tau}_Z={\varepsilon}_Z$, because
$U(y,{\eta})\cap Y=V(y,{\eta})\cap Y$ for $y\in Y$
and $U(z,{\eta})\cap Z=V(z,{\eta})\cap Z$ for $z\in Z$.
\eproof

\begin{definition}
Let $Y\subs \ncal$.
We say that a topological space $\<Y,{\tau}\>$
is a {\em standard refinement of $\<Y,{\varepsilon}\>$} provided that 
for each $y\in Y$ we can choose a  sequence of pairwise disjoint 
intervals with rational end points, $\ical^y=\{I^y_n:n\in {\omega}\}$,
which converges to $y$ such that  taking  $ J^y=\bigcup\ical^y$
the family 
$$
\bcal_y=\{U(y,{\eta})\setm  J^y:{\eta}>0\}.
$$
is a neighbourhood base of $y$ in ${\tau}$. 
\end{definition}

%It is easy to see that $X$ from theorem \ref{tm:zfc_example} is 
%a standard refinement.
%\end{comment}
\begin{theorem}\label{tm:ch_example}
If CH holds, then there is a  0-dimensional first countable 
standard refinement 
${\tau}$ of ${\varepsilon}$ on $\ncal$ such that
\begin{enumerate}\rlabel
\item $\RR(\<\ncal,{\tau}\>)={\omega}$,
\item $\nw(\<\ncal,{\tau}\>)=\cont$,
\item $\<\ncal,{\tau}\>$ has an irreducible base.
\end{enumerate}
\end{theorem}

\proof
First observe that for each $D\subs\ncal$  the set 
$\{x\in\overline{D}:x\notin\overline{(\infty,x)\cap D}\}$ 
is at most countable. 
Applying CH  and this observation  for each $y\in \ncal$  
we can choose a  sequence of pairwise disjoint
 intervals with rational endpoints, $\ical^y=\{I^y_n:n\in {\omega}\}$,
which is strictly increasing and converges to $y$,  such that taking
$ J^y=\bigcup\ical^y$ the assumptions  (A)--(B) below are satisfied:
\begin{enumerate}\Alabel
\item $\forall D\in\br\ncal;{\omega};$
$|\{y\in{\overline{D}}^{\varepsilon}:
y\notin\overline{D\cap J^y}^{\varepsilon}\}
|\le{\omega}$,
%\item $\forall D\in\br\ncal;{\omega};$
%$|\{y\in\overline{D}^{\varepsilon}:
%y\notin\overline{D\setm J^y}^{\varepsilon}\}
%|\le{\omega}$.
\end{enumerate}
To formulate property (B) we need the following notation:
for $D\in\br\ncal;{\omega};$ and $y\in \ncal$
write $\dy=\{d\in D:y\notin J^d\}$.
\begin{enumerate}\Alabel\addtocounter{enumi}{1}
\item $\forall D\in\br\ncal;{\omega};$
$|\{y\in\ncal:y\in\overline{\dy}^{\varepsilon}\land
%y\notin\overline{\dy\setm J^y}^{\varepsilon}\}
\dy\subs J^y
|\le{\omega}$.
\end{enumerate} 
Write $V(y,{\eta})=U(y,{\eta})\setm  J^y$ for ${\eta}>0$.
Let the neighbourhood base of $y$ in ${\tau}$ be 
$$
\bcal_y=\{V(y,{\eta}):{\eta}>0\}.
$$

Since 
\begin{equation}
\tag{\dag}\forall {\eta}>0\ \forall x\in V(y,{\eta})\setm \{y\}
\ \exists {\delta}>0\ U(x,{\delta})\subs V(y,{\eta})
\end{equation}
it follows that $\bcal=\{\bcal_y:y\in \ncal\}$ is a base of a topology. Since 
\begin{equation}
\tag{+}\forall {\eta}>0\ \forall{\delta}>0\ 
U(y,{\delta})\not\subs V(y,{\eta})
\end{equation}
it follows that the base $\bcal$ is irreducible.

It is not hard to see that (A)
implies that $\nw(\<\ncal,{\tau}\>)>{\omega}$.
Indeed, assume on the contrary that $\{M_m:m<{\omega}\}$ is a network.
Pick countable sets $K_m\subs M_m$ with 
$\overline{K_m}^{\varepsilon}=\overline{M_m}^{\varepsilon}$.
Then, by (A), there is $y\in\ncal$ such that for each $m\in {\omega}$
either $y\notin\overline{K_m}^{\varepsilon}=\overline{M_m}^{\varepsilon}$
or $y\in\overline{K_m\cap J^y}^{\varepsilon}
\subs\overline{M_m\cap J^y}^{\varepsilon}$.
Thus there is no $m\in {\omega}$ with $y\in M_m\subs \ncal\setm J^y$.

We will show that  (B)  implies  $\RR(\<\ncal,{\tau}\>)={\omega}$.
Assume on the contrary that $X$ is an uncountable weakly separated subspace
of $\<\ncal,{\tau}\>$. Since $\<\ncal,{\varepsilon}\>$ has countable weight,
we can assume that $x\in J^y$ or $y\in  J^x$ hold for each 
$x\ne y\in X$.

\renewcommand{\theclaim}{}
\begin{claim}
$\forall D\in \br X;{\omega};$
$|\{x\in X:x\in\overline{\dyi x;}^{\varepsilon}\}|\le {\omega}$.
\end{claim}   

\Proof{the claim}
If the above defined set is uncountable, then, by $(B)$,
there is $x\in X$ with $\dyi x;\not\subset J^x$. 
Let $d\in \dyi x;\setm J^x$ be arbitrary. 
Then $d\notin J^x$ and $x\notin J^d$ which contradicts 
our assumption on $X$. 
\Eproof

Using this claim, we can find an uncountable subset  
$Y=\{y_{\mu}:{\mu}<{\omega}_1\}$ of $X$ such that
$y_{\mu}\notin \overline {Y_{\mu}[y_{\mu}]}^{\varepsilon}$, where
$Y_{\mu}=\{y_{\nu}:{\nu}<{\mu}\}$.
So for each ${\mu}<{\omega}_1$ we have an interval $K_{\mu}$ with rational
endpoints such that  $y_{\mu}\in K_{\mu}$ and
for each ${\nu}<{\mu}$ if $y_{\nu}\in K_{\mu}$ then 
$y_{\nu}\notin Y_{\mu}[y_{\mu}]$, that is,
$y_{\mu}\in J^{y_{\nu}}$. Since there are only countable many 
intervals with rational endpoints, 
we can assume that $K_{\nu}=K$ for each ${\nu}<{\omega}_1$.
Since  $\ncal$ does not contain uncountable decreasing sequences, 
 there are ${\nu}<{\mu}<{\omega}_1$ with $y_{\nu}<y_{\mu}$.
But $J^y\subs (-\infty,y)$ by the construction, which contradicts to
$y_{\mu}\in J^{y_{\nu}}$.
So $\RR(\<\ncal,{\tau}\>)={\omega}$.
\eproof

%\begin{comment}
Let us remark that Todor{\v c}evi{\v c}, in \cite{To}, proved earlier that  
CH implies the existence of a  0-dimensional space $Y$ of size $\oo$ with
 $\w(Y)=\nw(Y)=\oo$  and ${\chi}(Y)=\RR(Y)={\omega}$.

The next theorem shows that some set-theoretic assumption
is necessary  to construct a standard refinement having the above described 
properties.  
To start with let us recall the Open Coloring Axiom (OCA)
(see \cite{To} and \cite{AS}).

\begin{ocaaxiom}
For each second countable  $T_3$ space  $X$ and 
open $H\subs \br X;2;$ either (i) or (ii) below holds:
\begin{enumerate}\rlabel
\item $X=\bigcup\limits_{n\in {\omega}} X_n$ where $X_n$ is $H$-independent,
\item $X$ contains an uncountable  $H$-complete subset. 
\end{enumerate}
\end{ocaaxiom}

\begin{theorem}\label{tm:oca}
$($OCA$)$ If $Y\subs\ncal$ and $\<Y,{\tau}\>$ is a standard refinement
of $\<Y,{\varepsilon}\>$ then either $\RR(Y)>{\omega}$ or
$\<Y,{\tau}\>$ is ${\sigma}$-second countable.
\end{theorem}

\proof
For each $y\in Y$ choose a  sequence of pairwise disjoint 
intervals with rational end points, $\ical^y=\{I^y_n:n\in {\omega}\}$,
which witnesses that $\<Y,{\tau}\>$ is a standard refinement.
Let  $ J^y=\bigcup\ical^y$.

Unfortunately the set 
{\em
$
E'=\{\<y,y'\>\in Y\times Y:y \in J_{y'}$ or $y'\in {J_y}\}$
} is not open in $Y\times Y$, so we need some extra work before applying OCA.

Fix an enumeration $\{K_k:k<{\omega}\}$ of the intervals with rational 
endpoints. For $y\in Y$ let us define the function
$f_y:{\omega}\to 2$ by taking $f_y(k)=1$ iff $K_k\subs J_y$.
Consider the second countable space 
$Z=\{\<y,f_y\>:y\in Y\}\subs Y\times \mbox{D(2)}^{\omega}$ and 
 define the set of edges $E$ on $Z$ as follows:
$$
\{\<y,f_y\>,\<y',f_{y'}\>\}\in E \Longleftrightarrow (y \in J_{y'} \mbox{ or }
y'\in {J_y}).
$$
It is easy to see that $E$ is open. So 
OCA implies that either there is an uncountable $E$-complete
$Z'\subs Z$ or $Z$ is the union of countable many 
$E$-independent subsets, $\{Z_n:n\in {\omega}\}$.

But by the definition of $E$, if $Z'$ is $E$-complete,
then $Y'=\{y\in Y:\<y,f_y\>\in Z'\}$ is weakly separated.
On the other hand, if $Z_n$ is $E$-independent, then 
${\tau}$ and ${\varepsilon}$ agree on $Y_n=\{y\in Y:\<y,f_y\>\in Z_n\}$.
\eproof

\begin{theorem}\label{tm:smtsd}
For each uncountable cardinal ${\kappa}$ there is a c.c.c poset $\pcalk$
of cardinality ${\kappa}$ such that in $V^{\pcalk}$ 
there is a   0-dimensional
first countable topological space $X=\<{\kappa},{\tau}\>$ and
there are c.c.c posets
$\qsd$ and $\qssc$   satisfying the following conditions: 
\begin{enumerate}\alabel
\item $V^{\pcalk}\models$ ``{\em $X$ has an irreducible base}'',
\item $V^{\pcalk*\qsd}\models$ ``{\em $X$ is ${\sigma}$-discrete}'', 
\item $V^{\pcalk*\qssc}\models$ ``{\em $X$ is ${\sigma}$-second countable}''. 
\end{enumerate}
So, in $V^{\pcalk}$,   $\w(X)={\kappa}$ by $($a$)$, $\nw$($X$)$={\kappa}$ by $($b$)$ and
$\RR(X^{\omega})={\omega}$ by (c).
\end{theorem}

\proof
We say that a quadruple $\anfg$ is in $\Pkoo$ provided
$(1)$--$(5)$ below hold:
\begin{enumerate}\arablabel
\item $A\in\br {\kappa};<{\omega};$,
\item $n\in {\omega}$,
\item $f$ and  $g$ are functions,
\item $f:A\times A\times n\to 2$,
\item $g:A\times n\times A\times n\to 3$,
\end{enumerate}

For $p\in \Pkoo$ we write $p=\anfgi p;$.
If $p, q\in \Pkoo$ we set
$p\le q$ iff   $f^p\supseteq f^q$ and $g^p\supseteq g^q$. 
If $p\in \Pkoo$, ${\alpha}\in \ap$, $i<\np$ set
$
U({\alpha},i)=\up({\alpha},i)=\{{\beta}\in \ap:\fp({\beta},{\alpha},i)=1\}
$.

A quadruple $\anfg\in \Pkoo$ is in $\Pko$ iff (i)--(iv) below are also 
satisfied:
\begin{enumerate}\rlabel
\item $\forall {\alpha}\in A$ $\forall i<n$ ${\alpha}\in U({\alpha},i)$,
\item $\forall {\alpha}\in A$ $\forall i<j<n$  
$U({\alpha},j)\subs U({\alpha},i)$,
\item $\forall {\alpha}\ne {\beta}\in A$ $\forall i,j<n$
\newline
\begin{tabular}{rcl}
%\newline\qquad
$U({\alpha},i)\subs U({\beta},j)$&iff& $g({\alpha},i,{\beta},j)=0$,\\
%\newline\qquad
$U({\alpha},i)\cap U({\beta},j)=\empt$&iff& $g({\alpha},i,{\beta},j)=1$.
\end{tabular}
\item $\forall {\alpha}\ne {\beta}\in A$ $\forall i,j<n$\newline
if ${\alpha}\in U({\beta},j)$ and ${\beta}\in U({\alpha},i)$
then $g({\alpha},i,{\beta},j)=2$.
\end{enumerate}
We claim that $\pcalk=\<\Pko,\le\>$ satisfies the requirements.

\begin{definition}
Assume that $p_i=\anfgi i;\in \Pkoo$ for $i\in 2$. We say that $p_0$ and $p_1$
are {\em twins} iff  $n_0=n_1$, $|A_0|=|A_1|$ and taking $n=n_0$ and
denoting by ${\sigma}$ the unique $<$-preserving bijection between $A_0$ and
$A_1$ we have
\begin{enumerate}\arablabel
\item ${\sigma}\rest {A_0\cap A_1}=\id_{A_0\cap A_1}$.
\item ${\sigma}$ is an isomorphism between $p_0$ and $p_1$, i.e.
$\forall {\alpha},{\beta}\in A_0$, $\forall i,j<n$
\begin{enumerate}\nolabel
\item $f_0({\alpha},{\beta},i)=f_1({\sigma}({\alpha}),{\sigma}({\beta}),i)$,
\item $g_0({\alpha},i,{\beta},j)=
g_1({\sigma}({\alpha}),i,{\sigma}({\beta}),j)$,
\end{enumerate}
\end{enumerate}
We say that ${\sigma}$ is the {\em twin function} of $p_0$ and $p_1$.
Define the {\em smashing function} $\sbar$ of $p_0$ and $p_1$ as follows:
 $\sbar={\sigma}\cup \id_{A_1}$.  
The function $\sstr$ defined by the formula 
$\sstr={\sigma}\cup {\sigma}^{-1}\rest {A_1}$ is called the
{\em exchange function} of $p_0$ and $p_1$.
\end{definition}

\begin{definition}
\label{def:eps_amalg}
Assume that $p_0$ and $p_1$ are twins and 
${\varepsilon}:A^{p_1}\setm A^{p_0}\to 2$. A common extension 
$q\in\Pko$ of $p_0$ and $p_1$  is called an
{\em ${\varepsilon}$-amalgamation} of the twins provided
$$
\forall {\alpha}\in A^{p_0}\triangle A^{p_1}\ 
f^q({\alpha},\sstr({\alpha}),i)={\varepsilon}(\sbar({\alpha})).
$$ 
\end{definition}

\begin{lemma}
\label{lm:twins}
If $p_0$,  $p_1\in\pcalk$ are twins and 
${\varepsilon}:A^{p_1}\setm A^{p_0}\to 2$, then $p_0$
and $p_1$ have an ${\varepsilon}$-amalgamation in $\Pko$.
\end{lemma}

\proof
Write $A=A_0\cup A_1$, $f^-=f_0\cup f_1$, $g^-=g_0\cup g_1$.
Let $B$ and $C$ be disjoint subsets of ${\kappa}\setm A$
of size $|A|$ and let ${\rho}:B\to A$ and ${\eta}:C\to A$ be 1-1.
Put $q=\<A\cup B\cup C,n,f,g\>$ where 
\begin{enumerate}\arablabel
\item $f^-\subs f$, $g^-\subs g$
\item $\forall {\alpha}\ne {\beta}\in A$ $\forall i,j<n$
$$f({\alpha},{\beta},i)= \left\{ 
\begin{array}{ll}
f_1(\sbar({\alpha}),\sbar({\beta}),i)&
\mbox{if $\sstr({\alpha})\ne {\beta}$},\\
{\varepsilon}(\sbar({\alpha}))&
\mbox{if $\sstr({\alpha})= {\beta}$,}
\end{array}
\right.
$$
and
$$
g({\alpha},i,{\beta},j)= \left\{ 
\begin{array}{ll}
g_1(\sbar({\alpha}),i,\sbar({\beta}),j)&
\mbox{if $\sstr({\alpha})\ne {\beta}$,}\\
2&\mbox{if $\sstr({\alpha})={\beta}$}.
\end{array}
\right.
$$
\item $\forall {\beta}\in B\cup C$ $\forall i<n$ $U^q({\beta},i)=\{{\beta}\}$,
\item $\forall {\alpha}\in A$ $\forall i,j<n$ $\forall {\beta}\in B$
$$
{\beta}\in U^q({\alpha},i)\mbox{\quad iff\quad}\exists l<n\, 
g({\rho}({\beta}),l,{\alpha},i)=0
$$
and
$$
g(\bj,\ai)= \left\{ 
\begin{array}{ll}
0&
\mbox{if ${\rho}({\beta})\in U^q({\alpha},i),$}\\
1&\mbox{if ${\rho}({\beta})\notin U^q({\alpha},i)$}.
\end{array}
\right.
$$
\item $\forall {\alpha}\in A$ $\forall i,l<n$ $\forall {\gamma}\in C$
$$
{\gamma}\in U^q({\alpha},i)\mbox{\quad iff\quad} 
\sbar({\eta}({\gamma}))\in U_1(\sbar({\alpha}),i).
$$
and
$$
g(\cl,\ai)= \left\{ 
\begin{array}{ll}
0&
\mbox{if ${\eta}({\gamma})\in U^q({\alpha},i),$}\\
1&\mbox{if ${\eta}({\gamma})\notin U^q({\alpha},i)$}.
\end{array}
\right.
$$
\end{enumerate}

Let us remark that $(4)$ and $(5)$ contain no
circularity  because $(1)$ and $(2)$ 
define $g$ on $A$.
Obviously $q\in \Pkoo$ and $q\le p_0,p_1$, so we have to show that $q\in \Pko$.
(i) is straightforward. Before checking (ii)--(iv) we need some 
preparation.
If ${\alpha},{\beta}\in A$ and $i,j<n$ write $\<\ai\>\tle\<\bj\>$
iff $g(\ai,\bj)=0$.
\begin{claim}
\label{claim:transitive}
The relation $\tle$ is transitive on $A\times n$.
\end{claim}

\Proof{claim}
Assume that $\<\ai\>\tle\<\bj\>\tle\<\cl\>$.
Then, by $(2)$, $\<\sai\>\tle\<\sbj\>\tle\<\scl\>$, so
$\<\sai\>\tle\<\scl\>$. Thus $\<\ai\>\tle\<\cl\>$ provided $\ssa\ne{\gamma}$.
 Assume  that $\ssa={\gamma}$.  Then
$U_1(\sba,i)\subs U_1(\sbb,j)\subs U_1(\sba,l)$ and so
$\sba\in U_1(\sbb,j)$ and $\sbb\in U_1(\sba,l)$. Thus
$g_1(\sbb,j,\sba,l)=2$ because $p_1$ satisfies (iv).
This contradiction proves that  $\ssa\ne{\gamma}$.
\Eproof

\begin{claim}
\label{claim:zero}
$\forall {\alpha}\ne {\beta}\in A$ $\forall i,j<n$
if $g(\ai,\bj)=0$ then $U^q(\ai)\subs U^q(\bj)$.
\end{claim}

\Proof{claim}
We have $g_1(\sba,i,\sbb,j)=0$. 
Thus 
\begin{equation}\tag{\dag}U_1(\sba,i)\subs U_1(\sbb,j).
\end{equation}
Let ${\gamma}\in A\cap U^q({\alpha},i)$. 
If $\ssc\ne {\beta}$ then (\dag) implies ${\gamma}\in U^q({\beta},j)$
by (2).
Assume now that $\ssc={\beta}$.
Then $\sbb=\sbc\in U_1(\sba,i)$ and, on the other hand,
$\sba\in U_1(\sbb,j)$  by (\dag).
Thus $g_1(\sba,i,\sbb,j)=2$ because $p_1$ satisfies (iv).
Contradiction, $\ssc\ne{\beta}$.
So we have shown 
$U^q(\ai)\cap A\subs U^q(\bj)\cap A$. 
Next we can see that 
$U^q(\ai)\cap B\subs U^q(\bj)\cap B$ by claim \ref{claim:transitive} 
and by (4).
Finally let ${\gamma}\in U^q({\alpha},i)\cap C$.
Then $\sbar({\eta}({\gamma}))\in U_1(\sba,i)$, so
$\sbar({\eta}({\gamma}))\in U_1(\sbb,j)$. Thus
${\gamma}\in U^q({\beta},j)$ by (5). 
\Eproof

\begin{claim}
\label{claim:one}
$\forall {\alpha}\ne {\beta}\in A$ $\forall i,j<n$
if $g(\ai,\bj)=1$ then $U^q(\ai)\cap U^q(\bj)=\empt$.
\end{claim}

\Proof{claim}
Since $g(\ai,\bj)=1$ we have $\ssa\ne{\beta}$, so 
$U_1(\sba,i)\cap U_1(\sbb,j)=\empt$ implies
 $U^q(\ai)\cap U^q(\bj)\cap A=\empt$.
Assume now that ${\gamma}\in U^q(\ai)\cap U^q(\bj)\cap B$.
Then ${\rho}({\gamma})\in U^q(\ai)\cap U^q(\bj)\cap A$, contradiction.
Finally assume that ${\gamma}\in U^q(\ai)\cap U^q(\bj)\cap C$.
Then $\sbar({\eta}({\gamma}))\in U_1(\sba,i)\cap U_1(\sbb,j)=\empt$, 
which is impossible.
\Eproof

\begin{claim}
\label{claim:two}
$\forall {\alpha}\ne {\beta}\in A\cup B\cup C$ $\forall i,j<n$
if $g(\ai,\bj)=2$ then 
$\empt\ne U^q(\ai)\setm U^q(\bj)\ne U^q(\ai)$.
\end{claim}

\Proof{claim}
The assumption implies ${\alpha},{\beta}\in A$. If
${\alpha}\ne \sstr({\beta})$ then $g(\sba,i,\sbb,j)=2$, which
implies the statement.
So we can assume that ${\alpha}=\sstr({\beta})$. Let 
${\gamma}={\rho}^{-1}({\alpha})$. Then ${\gamma}\in U^q(\ai)$
and, on the other hand, ${\gamma}\notin U^q(\bj)$,
because $g({\alpha},l,\bj)=2$ for each $l<n$ by $(2)$.
So ${\gamma}\in U^q(\ai)\setm U^q(\bj)$.
Let ${\delta}={\eta}^{-1}({\alpha})$. 
Then $\sbar ({\eta}({\delta}))=\sbb=\sba\in U_1(\sba,i)\cap U_1(\sbb,j)$,
so ${\delta}\in U^q(\ai)\cap U^q(\bj)$.
\Eproof

So we have proved (iii) for $q$, which implies (ii).
To check (iv) assume that ${\alpha}\ne {\beta}\in A\cup B\cup C$,
$i,j<n$ with  ${\alpha}\in U^q(\bj)$ and ${\beta}\in U^q(\ai)$.
Then ${\alpha},{\beta}\in A$. If $\ssa\ne {\beta}$, then
$\sba\in U_1(\sbb,j)$ and $\sbb\in U_1(\sba,i)$ implies
$g_1(\sba,i,\sbb,j)=2$. Hence $g(\ai,\bj)=2$.
If $\ssa= {\beta}$, then $g(\ai,\bj)=2$ by definition.
The lemma is proved.
\eproof

The previous lemma implies that $\pcalk$ satisfies c.c.c 
because among uncountably many elements of $\pcalk$ there are always two twins.
Let $\gcal$ be the $\pcalk$ generic filter and
let $F=\bigcup\{f^p:p\in\gcal\}$.
For each ${\alpha}<{\kappa}$ and $n\in{\omega}$ let
$V({\alpha},i)=\{{\beta}<{\kappa}:F({\beta},{\alpha},i)=1\}$.
Put $\bcal_{\alpha}=\{V({\alpha},i):i<{\kappa}\}$ and
$\bcal=\bigcup\{\bcal_{\alpha}:{\alpha}<{\kappa}\}$.
By standard density arguments we can see that
$\bcal$ is  base of a first countable topological space 
$X=\<{\kappa},{\tau}\>$. Since $\pcalk$ satisfies (iv), 
$\bigcup\{\bcal_{\alpha}:{\alpha}<{\kappa}\}$ is 
an irreducible decomposition of $\bcal$.

We are now ready to define the posets $\qsd$ and $\qssc$
in $V^{\pcalk}$.

A triple $\bde$ is in  $\qsd$ iff
\begin{enumerate}\alabel
\item $B\in\br{\kappa};<{\omega};$,
\item $d:B\to {\omega}$,
\item $e:B\to {\omega}$,
\item $\forall {\alpha}\ne {\beta}\in B$ if $d({\alpha})=d({\beta})$
then ${\alpha}\notin V({\beta},e({\beta}))$.
\end{enumerate}

A quadruple $\bdne$ is in $\qssc$ iff
\begin{enumerate}\Alabel
\item $B\in\br{\kappa};<{\omega};$,
\item $d:B\to {\omega}$,
\item $m\in{\omega}$,
\item $e:B\times m\to {\omega}$,
\item $\forall {\alpha},{\beta},{\gamma}\in B$ 
$\forall i,j<m$ if $d({\alpha})=d({\beta})=d({\gamma})$ and
$e({\alpha},i)=e({\beta},j)$ then
 ${\gamma}\in V({\alpha},i)$ iff ${\gamma}\in V({\beta},j)$.
\end{enumerate}
The orderings on $\qsd$ and  $\qssc$ are defined in the straightforward 
way. If $q$ and $r$ are compatible elements of $Q_i$, then denote by
$q\land r$ their greatest lower bound in $Q_i$.
\begin{lemma}
$\pcalk*\qsd$ satisfies c.c.c.
\end{lemma}

\proof
Let $\<\<p_{\nu},q_{\nu}\>:{\nu}<\oo\>\subs \pcalk*\qsd$.
We can assume that $p_{\nu}$ decides $q_{\nu}$. Write
$p_{\nu}=\anfgi {\nu};$ and $q_{\nu}=\bdei {\nu};$.
By standard density arguments we can assume that
$A^{\nu}\supset B^{\nu}$. Applying standard $\Delta$-system and counting
arguments  we can find 
${\nu}<{\mu}<\oo$ such that $p_{\nu}$ and $p_{\mu}$ are twins 
and denoting by ${\sigma}$ the twin function of $p_{\nu}$ and
$p_{\mu}$ we have $d^{\nu}({\alpha})=d^{\mu}({\sigma}({\alpha}))$
and $e^{\nu}({\alpha})=e^{\mu}({\sigma}({\alpha}))$
for each ${\alpha}\in B^{\nu}$.

Define the function ${\varepsilon}^0:A^{\mu}\setm A^{\nu}\to 2$ by the
equation ${\varepsilon}^0({\alpha})=0$. By lemma \ref{lm:twins}
the conditions $p^{\nu}$ and $p^{\mu}$ have an 
${\varepsilon}^0$-amalgamation $p$.
We claim that 
$$
p\force q_{\nu}\land q_{\mu}\in \qsd.
$$
It is enough to check that if ${\alpha}\in B^{\nu}$,
${\beta}\in B^{\mu}$, $d^{\nu}({\alpha})=d^{\mu}({\beta})$, then
$p\force {\alpha}\notin V({\beta},e^{\mu}({\beta}))$.
If $\ssa\ne {\beta}$, then (2) implies this.
If $\ssa={\beta}$, then $p\force {\alpha}\notin V({\beta},e^{\mu}({\beta}))$
because $p$ is an ${\varepsilon}^0$-amalgamation.
\eproof

\begin{lemma}
$\pcalk*\qssc$ satisfies c.c.c.
\end{lemma}

\proof
Let $\<\<p_{\nu},q_{\nu}\>:{\nu}<\oo\>\subs \pcalk*\qssc$.
We can assume that $p_{\nu}$ decides $q_{\nu}$. Write
$p_{\nu}=\anfgi {\nu};$ and $q_{\nu}=\bdnei {\nu};$.
By standard density arguments we can assume that
$A^{\nu}\supset B^{\nu}$. Applying standard arguments we can find 
${\nu}<{\mu}<\oo$ such that $p_{\nu}$ and $p_{\mu}$ are twins 
and denoting by ${\sigma}$ the twin function of $p_{\nu}$ and
$p_{\mu}$ we have $m^{\nu}=m^{\mu}$ and
$d^{\nu}({\alpha})=d^{\mu}({\sigma}({\alpha}))$ and
$e^{\nu}({\alpha},i)=e^{\mu}({\sigma}({\alpha}),i)$
for each ${\alpha}\in A^{\nu}$ and $i<m^{\mu}$.

Define the function ${\varepsilon}^1:A^{\mu}\setm A^{\nu}\to 2$ by the
equation ${\varepsilon}^1({\alpha})=1$. By lemma \ref{lm:twins}
the conditions $p^{\nu}$ and $p^{\mu}$ have an 
${\varepsilon}^1$-amalgamation $p$.
We claim that 
$$
p\force q^{\nu}\land q^{\mu}\in \qssc.
$$
Set $d=d^{\nu}\cup d^{\mu}$ and $e=e^{\nu}\cup e^{\mu}$.
Let ${\alpha},{\beta},{\gamma}\in B^{\nu}\cup B^{\mu}$
with $d({\alpha})=d({\beta})=d({\gamma})$ and
$e({\alpha},i)=e({\beta},j)$.
We have to show that
\begin{equation}
\tag{$*$}p\force\mbox{ \quad } 
{\gamma}\in V(\ai)\mbox{\quad iff\quad} {\gamma}\in V(\bj).
\end{equation}
We know that $p\force$ {\em ``$\sbc\in V(\sba,i)$ iff $\sbc\in V(\sbb,j)$''},
so we can assume that 
$p\force$ {\em ``$\sbc\in V(\sba,i)$ and $\sbc\in V(\sbb,j)$''},.
So if $\ssc$ is different from ${\alpha}$ and ${\beta}$, then we are done.
Assume finally that $\ssc={\alpha}$. Then 
$p\force$ {\em ``$\sbc\in V(\sba,i)$''}, so
$p\force$ {\em ``$\sbc\in V(\sbb,j)$''}, thus 
$p\force$ {\em ``${\gamma}\in V(\bj)$''} by (2).
But $p$ is an ${\varepsilon}^1$-amalgamation, so 
$p\force $ {\em ``$\ssa\in V({\alpha},i)$''}, i.e. 
$p\force$ {\em ``${\gamma}\in V(\ai)$''}.
\eproof

\begin{lemma}
$$
V^{\pcalk*\qsd}\models\mbox{``$X$ is ${\sigma}$-discrete''}
$$
\end{lemma}

\proof
Let $\hcal$ be the $\qsd$-generic filter over $V^{\pcalk}$.
Set $d=\bigcup\{d^q:q\in \hcal\}$ and
$e=\bigcup\{e^q:q\in \hcal\}$. By standard density arguments 
the domains of the functions  $d$ and $e$ are ${\kappa}$.
We have $V(x,e(x))\cap d^{-1}\{d(x)\}=\{x\}$ by (d), so
 $d^{-1}(n)$ is discrete for each $n\in {\omega}$.
\eproof

\begin{lemma}
$$
V^{\pcalk*\qssc}\models\mbox{``$X$ is ${\sigma}$-second countable''}
$$
\end{lemma}

\proof
Let $\hcal$ be the $\qssc$-generic filter over $V^{\pcalk}$.
Take $d=\bigcup\{d^q:q\in \hcal\}$ and
$e=\bigcup\{e^q:q\in \hcal\}$. 
By standard density arguments 
$\dom(d)={\kappa}$ and $\dom(e)={\kappa}\times {\omega}$.
Fix $n\in {\omega}$ and let $X_n=d^{-1}\{n\}$.
We claim that $\w(X_n)={\omega}$. Indeed,
$\bigcup\{V({\alpha},i):{\alpha}\in X_n,i<{\omega}\}$ is a base
of $X_n$ and by (E), if ${\alpha},{\beta}\in X_n,$
$i,j<{\omega}$ and $e({\alpha},i)=e({\beta},j)$ then
$V({\alpha},i)\cap X_n=V({\beta},j)\cap X_n$. 
\eproof
This completes the proof of the theorem.
\eproof

We have shown in \cite{JSS2} that every first countable $T_2$ space
satisfying $\RR(X^{\omega})={\omega}$ becomes ${\sigma}$-second countable
in a suitable c.c.c extension. Thus (c) of theorem \ref{tm:smtsd}
may be considered as the natural way to insure $\RR(X^{\omega})={\omega}$.

%\input{mbigwref}

%\end{comment}

\end{document}